\documentclass[12pt,reqno]{amsart}

\usepackage{amssymb}
\usepackage{amsmath}
\usepackage{amsthm}
\usepackage{mathrsfs}
\usepackage{bm}

\numberwithin{equation}{section}

\theoremstyle{plain}
\newtheorem{thm}{Theorem}[section]
\newtheorem{prop}[thm]{Proposition}
\newtheorem{lem}[thm]{Lemma}

\theoremstyle{definition}

\newcommand{\ppp}{\partial}
\newcommand{\A}{\mathcal{A}}
\newcommand{\B}{\mathcal{B}}
\newcommand{\R}{\mathbb{R}}
 
\newcommand{\N}{\mathbb{N}}

\providecommand{\abs}[1]{\left\lvert#1\right\rvert}
\providecommand{\norm}[1]{\left\lVert#1\right\rVert}

\allowdisplaybreaks

\begin{document}

\title[Determination of coefficients in fractional diffusion equations]
{Determination of time dependent factors of coefficients in 
fractional diffusion equations}

\author[Kenichi Fujishiro]{Kenichi Fujishiro${}^{1)}$}
\author[Yavar Kian]{Yavar Kian${}^{2)}$}

\thanks{1) Department of Mathematical Sciences, University of 
Tokyo, Komaba, Meguro, Tokyo 153, Japan}
\thanks{2) Aix-Marseille Universit\'e, CNRS, CPT UMR 7332, 13288 Marseille, 
France and Universit\'e de Toulon, CNRS, CPT UMR 7332 83957 La Garde, France}



\email{1) kenichi@ms.u-tokyo.ac.jp}
\email{2) yavar.kian@univ-amu.fr}

\date{}

\begin{abstract}
  In the present paper, we consider initial-boundary value problems for 
partial differential equations with time-fractional derivatives which evolve 
in $Q=\Omega\times(0,T)$ where $\Omega$ is a bounded domain of $\R^d$ and 
$T>0$.
  We study the stability of the inverse problem of determining the 
time-dependent parameter in a source term or a coefficient of zero-th order 
term from observations of the solution at a point $x_0\in\overline{\Omega}$ 
for all $t\in(0,T)$.
\end{abstract}

\maketitle

\section{Introduction}

  Let $\Omega$ be a bounded domain of $\R^d$, $d=1,2,3$, with $\mathcal C^2$ 
boundary $\ppp\Omega$.
  We set $\Sigma=\ppp\Omega\times(0,T)$ and $Q=\Omega\times(0,T)$.
  We consider the following two initial-boundary value problem (IBVP in 
short) for the fractional diffusion equation
\begin{equation}\label{eq1}
\left\{\begin{aligned}
	\ppp_t^\alpha u(x,t)+\A u(x,t)=f(t)R(x,t), \quad &(x,t)\in Q, \\
	\B_{\sigma} u(x,t)=0,	\quad &(x,t)\in\Sigma,\\
	u(x,0)=0,	\quad &x\in\Omega
\end{aligned}\right.
\end{equation}
  and
\begin{equation}\label{eq2}
\left\{\begin{aligned}
	\ppp_t^\alpha v(x,t)+\A v(x,t)+f(t)q(x,t)v(x,t)=0, \quad&(x,t)\in Q,\\
	\B_{\sigma} v(x,t)=0,	\quad&(x,t)\in\Sigma, \\
	v(x,0)=v_0(x),	\quad &x\in\Omega
\end{aligned}\right.
\end{equation}
  with $0<\alpha<1$.
  Here we denote by $\ppp_t^\alpha$ the Caputo fractional derivative with 
respect to $t$:
\[
	\ppp_t^\alpha g(t)
	:=\frac{1}{\Gamma(1-\alpha)}\int_0^t(t-s)^{-\alpha}\frac{dg}{ds}(s)ds.
\]
  The differentiual operator $\A$ is defined by
\[
	\A u(x,t):=
	-\sum_{i,j=1}^d\frac{\ppp}{\ppp x_i}
		\left(a_{ij}(x)\frac{\ppp u}{\ppp x_j}(x,t)\right),
\]
  where the coefficients satisfy
\[
	a_{ij}=a_{ji},\quad 1\le i,j\le d,
	\quad\mbox{and}\quad
	\sum_{i,j=1}^d a_{ij}(x)\xi_i\xi_j\ge\mu|\xi|^2,
		\quad x\in\overline{\Omega},\ \xi\in\R^d
\]
  for some $\mu>0$.
  Moreover $\B_{\sigma}$ is defined as
\[
	\B_{\sigma} u(x)=(1-\sigma(x))u(x)+\sigma(x)\ppp_{\nu_A}u(x),
	\quad x\in\ppp\Omega,
\]
  where 
\[
	\ppp_{\nu_A}u(x)
	=\sum_{i,j=1}^d a_{ij}(x)\frac{\ppp u}{\ppp x_i}\nu_j(x)
\]
  and $\nu=(\nu_1,\dots,\nu_d)$ is the outward normal unit vector to 
$\ppp\Omega$.
  Here $\sigma$ is a $\mathcal C^2$ function on $\ppp\Omega$ satisfying
\[
	0\le\sigma(x)\le1,\quad x\in\ppp\Omega.
\]
  For the regularity of $a_{ij}$, we assume
\begin{equation*}\label{aij}
\begin{cases}
	a_{ij}\in\mathcal C^1(\overline{\Omega}) &\mbox{if}\ \sigma\equiv0, \\
	a_{ij}\in\mathcal C^2(\overline{\Omega}) &\mbox{if}\ \sigma\not\equiv0.
\end{cases}
\end{equation*}
  Note that the regularity for $a_{ij}$ depends on whether $\sigma\equiv0$ or 
not, which is due to condition \eqref{Aeq2} in the next section.

  In the present paper, we consider the inverse problem which consists of 
determining the function $\{f(t)\}_{t\in(0,T)}$ in \eqref{eq1} and \eqref{eq2} 
from the observation of the solution at a point $x_0\in\overline{\Omega}$ for 
all $t\in(0,T)$.

  The partial differential equations with time fractional derivatives such as 
\eqref{eq1} and \eqref{eq2} are proposed as new models describing the 
anomalous diffusion phenomena.
  Adams and Gelhar \cite{AG} pointed out that the field data in a highly 
heterogeneous aquifer cannot be described well by the classical advection 
diffusion equation.
  Hatano and Hatano \cite{HA} applied the continuous-time random walk (CTRW) 
as a microscopic model of the diffusion of ions in heterogeneous media.
  From the CTRW model, one can derive a fractional diffusion equation as a 
macroscopic model (see e.g., Metzler and Klafter \cite{MK} and Roman and 
Alemany \cite{RA}).
  In particular, the fractional diffusion equation can be used as a model for 
the diffusion of contaminants in a soil.
  Therefore the inverse problem considered in this paper means the 
determination of the time evolution of pollution source in \eqref{eq1} and 
reaction rate of pollutants in \eqref{eq2} respectively.
  In this paper, we consider such problems assuming the boundedness of the 
time-dependent parameter $\{f(t)\}_{t\in(0,T)}$ (see \eqref{ft}).

  As monographs of fractional calculus, there are books such as Podlubny 
\cite{P} and Samko, Kilbas and Marichev \cite{SM} for example.
  As for mathematical works concerned with partial differential equations 
with time fractional derivatives, we can refer to Agarwal \cite{AW}, Gejji and 
Jafari \cite{GJ}, Gorenflo and Mainardi \cite{GM}, Luchko \cite{L} and 
references therein.

\bigskip
  The remainder of this paper is composed of four sections.
  In Section 2, we state our main results.
  In Section 3, we study the forward problem for the IBVPs \eqref{eq1} and 
\eqref{eq2} and prove the unique existence and regularity of the solutions.
  In Sections 4 and 5, we complete the proof of our main results---Theorems 
\ref{t1} and \ref{t2} respectively.

\section{Main results}

  By $L^2(\Omega)$, we denote the usual $L^2$-space equipped with the inner 
product $(\cdot,\cdot)$ and the norm $\|\cdot\|:=\|\cdot\|_{L^2(\Omega)}$.
  Moreover $H^s(\Omega)$, $s\in\R$, and $W^{m,p}(\Omega)$, $m\in\N$, 
$1\le p\le\infty$, are the Sobolev spaces (see Adams \cite{AD} for example).

  For the time dependent parameter $\{f(t)\}_{t\in(0,T)}$, we always assume 
\begin{equation}\label{ft}
	f\in L^{\infty}(0,T).
\end{equation}
  For other given functions in \eqref{eq1}, we suppose
\begin{equation}\label{Aeq1}
	R\in L^p(0,T;H^2(\Omega)),\ \frac{8}{\alpha}<p\le\infty
	\quad\mbox{and}\quad
	\B_{\sigma}R=0\quad\mbox{on}\ \Sigma.
\end{equation}
  On the other hand, in the IBVP \eqref{eq2}, we suppose
\begin{equation}\label{Aeq2}
\begin{cases}
	q\in L^{\infty}(0,T;H^{2}(\Omega))\quad
	(\mbox{and $\ppp_{\nu}q=0$ on $\Sigma$ if $\sigma\not\equiv0$}),\\
	v_0\in H^4(\Omega)\quad\mbox{and}\quad 
	\B_{\sigma}v_0=\B_{\sigma}(\A v_0)=0\quad\mbox{on}\ \ppp\Omega.
\end{cases}
\end{equation}
  Assuming these conditions, we prove in Section 3 that the IBVPs 
\eqref{eq1} and \eqref{eq2} admit unique solutions 
$u,v\in\mathcal C([0,T]; H^{2}(\Omega))$ with 
$\ppp_t^\alpha u\in L^p(0,T;H^s(\Omega))$ and 
$\ppp_t^\alpha v\in L^{\infty}(0,T;H^s(\Omega))$ for some $s>d/2$.
  Therefore, using the Sobolev embedding theorem (see Theorem 9.8 in Chapter 1 
of \cite{LM} for example), for any $x_0\in\overline{\Omega}$, we see that 
\[
	\ppp_t^\alpha u(x_0,\cdot)\in L^p(0,T)
	\quad\mbox{and}\quad
	\ppp_t^\alpha v(x_0,\cdot)\in L^{\infty}(0,T).
\]
  Then our main results can be stated as follows;

\bigskip
\begin{thm}\label{t1}
  Let condition \eqref{Aeq1} be fulfilled and $u_i$ be the solution of 
\eqref{eq1} for $f=f_i\in L^{\infty}(0,T)$, ($i=1,2$).
  We assume that there exist $x_0\in\overline{\Omega}$ and $\delta>0$ such that
\begin{equation}\label{t1a}
	|R(x_0,t)|\geq \delta,\quad\mbox{a.e.}\ t\in(0,T).
\end{equation}
  Then there exists a constant $C>0$ depending on $p$, $T$, $\Omega$, $\delta$ 
and $\|R\|_{L^p(0,T;H^2(\Omega))}$ such that
\begin{align}
&	\|f_1-f_2\|_{L^p(0,T)}
	\le C\|\ppp_t^\alpha u_1(x_0,\cdot)
		-\ppp_t^\alpha u_2(x_0,\cdot)\|_{L^p(0,T)}, \label{t1b} \\
&	\|\ppp_t^\alpha u_1(x_0,\cdot)
		-\ppp_t^\alpha u_2(x_0,\cdot)\|_{L^p(0,T)}
\le	C\|f_1-f_2\|_{L^{\infty}(0,T)}. \label{t1b'}
\end{align}
  In particular, if we take $p=\infty$ in \eqref{Aeq1}, then
\begin{align*}
	C^{-1}\|\ppp_t^\alpha u_1(x_0,\cdot)
		-\ppp_t^\alpha u_2(x_0,\cdot)\|_{L^{\infty}(0,T)}
&\le	\|f_1-f_2\|_{L^{\infty}(0,T)} \\
&\le 	C\|\ppp_t^\alpha u_1(x_0,\cdot)
		-\ppp_t^\alpha u_2(x_0,\cdot)\|_{L^{\infty}(0,T)}.
\end{align*}
\end{thm}

\bigskip
\begin{thm}\label{t2}
  Let condition \eqref{Aeq2} be fulfilled and $v_i$ be the 
solution of \eqref{eq2} for $f=f_i\in L^{\infty}(0,T)$ with 
$\|f_i\|_{ L^{\infty}(0,T)}\le M$ ($i=1,2$).
  We assume that there exist $x_0\in\overline{\Omega}$ and $\delta>0$ such that
\begin{equation}\label{t2a}
	|q(x_0,t)v_2(x_0,t)|\geq \delta,\quad\mbox{a.e.}\ t\in(0,T).
\end{equation}
  Then there exists a constant $C>0$ depending on $M$, $T$, $\Omega$, 
$\delta$ and $\|q\|_{L^{\infty}(0,T;H^{2}(\Omega))}$ such that
\begin{align}\label{t2b}
	C^{-1}\|\ppp_t^\alpha v_1(x_0,\cdot)
		-\ppp_t^\alpha v_2(x_0,\cdot)\|_{L^{\infty}(0,T)}
	&\le\|f_1-f_2\|_{L^{\infty}(0,T)} \nonumber \\
	&\le C\|\ppp_t^\alpha v_1(x_0,\cdot)
		-\ppp_t^\alpha v_2(x_0,\cdot)\|_{L^{\infty}(0,T)}.
\end{align}
\end{thm}

\bigskip
  In Theorem 4.4 of Sakamoto and Yamamoto \cite{SY}, a similar problem to 
Theorem \ref{t1} is considered, but our result is more applicable in the point 
of view that the factor $R(x,t)$ is also allowed to depend on $t$.
  Moreover, we may assume less regularity for $R$ in Theorem \ref{t1}.
  The arguments of Theorem \ref{t2} can also be applied to parabolic equations 
in order to consider the result of Theorem 1.1 in \cite{CK} with observations 
of the solution at a point $x_0\in\Omega$ when $\Omega\subset\R^d$, $d=1,2,3$.

  For such inverse problems with $\alpha=1$, we can also refer to Section 1.5 
of Prilepko, Orlovsky and Vasin \cite{PR}, Cannon and Esteva \cite{CE} and 
Saitoh, Tuan and Yamamoto \cite{STY1,STY2}, for example.
  In our main results, we assume conditions \eqref{t1a} and \eqref{t2a}, which 
means that the observation point cannot be far from the source.
  On the other hand, in \cite{CE} and \cite{STY1,STY2}, the 
determination of time dependent factor in the source term is considered 
without assuming such conditions and the logarithmic type and H\"older type 
estimates are proved respectively.
  However, the results for fractional diffusion equations without these 
conditions have not been obtained yet.
  Here we restrict ourselves to the case with assumptions \eqref{t1a} and 
\eqref{t2a}, and show the Lipschitz type stability.

  Let us remark that the results of this paper can be extended to the case 
$d\ge4$.
  For this purpose additional conditions such as more regularity for $a_{ij}$ 
and $\ppp\Omega$ are required.
  In order to avoid technical difficulties, we only treat the case $d\le3$.

\section{Forward problem}

  This section is devoted to the proof of unique existence and regularity of 
the solution of the IBVPs \eqref{eq1} and \eqref{eq2}.

\bigskip
\begin{prop}\label{p1}
  Let conditions \eqref{ft} and \eqref{Aeq1} be fulfilled.
  Then the IBVP \eqref{eq1} admits a unique solution 
$u\in\mathcal C([0,T];H^{2}(\Omega))$  satisfying
\[
	\A u\in\mathcal C([0,T];H^{2\gamma}(\Omega))
	\quad\mbox{and}\quad
	\ppp_t^\alpha u\in L^p(0,T;H^{2\gamma}(\Omega))
\]
  for all $0\le\gamma<1-1/(p\alpha)$.
  Moreover we have 
\begin{align}\label{p1a}
	\|\A u\|_{\mathcal C([0,T];H^{2\gamma}(\Omega))}
		+\|\ppp_t^\alpha u\|_{L^p(0,T;H^{2\gamma}(\Omega))}
\le 	C\|fR\|_{L^p(0,T;H^2(\Omega))}.
\end{align}
  with $C>0$ depending on $\Omega$, $T$ and $\gamma$
\end{prop}

\bigskip
\begin{prop}\label{p2}
  Let conditions \eqref{ft} and \eqref{Aeq2} be 
fulfilled.
  Then the IBVP \eqref{eq2} admits a unique solution 
$v \in\mathcal C([0,T];H^{2}(\Omega))$  satisfying
\[
	\A v\in\mathcal C([0,T];H^{2\gamma}(\Omega))
	\quad\mbox{and}\quad
	\ppp_t^\alpha v\in L^{\infty}(0,T;H^{2\gamma}(\Omega))
\]
  for all $0\le\gamma<1$.
  Moreover, we have
\begin{equation}\label{p2a}
	\|\A v\|_{\mathcal C([0,T];H^{2\gamma}(\Omega))}
	+\|\ppp_t^\alpha v\|_{L^{\infty}(0,T;H^{2\gamma}(\Omega))}
	\leq C\|v_0\|_{H^4(\Omega)}
\end{equation}
  with $C$ depending on $\Omega$, $T$, $\|f\|_{L^{\infty}(0,T)}$, 
$\|q\|_{L^{\infty}(0,T;H^{2}(\Omega))}$ and $\gamma$.
\end{prop}

\bigskip
  If all coefficients are independent of time variable $t$, then we can apply 
eigenfunction expansion and the problems can be reduced to ordinary 
differential equations of fractional order (e.g. \cite{SY}).
  However, since we consider the determination of the time dependent factor of 
coefficients, we apply fixed point theorem to show the unique existence of the 
solutions to \eqref{eq1} and \eqref{eq2} as in Beckers and Yamamoto \cite{BY}.

  In order to prove these results, we consider the IBVPs with more general 
data in the next subsections.

\subsection{Intermediate results}
  We introduce the following IBVPs
\begin{equation}\label{eq3}
\left\{\begin{aligned}
	\ppp_t^\alpha u(x,t)+\A u(x,t)=F(x,t),	\quad &(x,t)\in Q,\\
	\B_{\sigma} u(x,t)=0,			\quad &(x,t)\in\Sigma, \\
	u(x,0)=0,				\quad &x\in\Omega,
\end{aligned}\right.
\end{equation}
\begin{equation}\label{eq4}
\left\{\begin{aligned}
	\ppp_t^\alpha v(x,t)+\A v(x,t)+p(x,t)v(x,t)=F(x,t),\quad &(x,t)\in Q,\\
	\B_{\sigma} v(x,t)=0,			\quad &(x,t)\in\Sigma, \\
	v(x,0)=0,			\quad &x\in\Omega,
\end{aligned}\right.
\end{equation}
  and
\begin{equation}\label{eq5}
\left\{\begin{aligned}
	\ppp_t^\alpha v(x,t)+\A v(x,t)+p(x,t)v(x,t)=0,	\quad &(x,t)\in Q,\\
	\B_{\sigma} v(x,t)=0,			\quad &(x,t)\in\Sigma, \\
	v(x,0)=v_0(x),			\quad &x\in\Omega.
\end{aligned}\right.
\end{equation}
  We also consider the following conditions
\begin{equation}\label{Aeq1(bis)}
	F\in L^p(0,T;H^2(\Omega)),\ \frac{8}{\alpha}<p\le\infty
	\quad\mbox{and}\quad
	\B_{\sigma}F=0\quad\mbox{on}\ \Sigma
\end{equation}
  and
\begin{equation}\label{Aeq2(bis)}
\begin{cases}
	1)\ p\in L^{\infty}(0,T;H^{2}(\Omega)) \quad
	(\mbox{and $\ppp_{\nu}p=0$ on $\Sigma$ if $\sigma\not\equiv0$}), \\
	2)\ v_0\in H^4(\Omega)\quad\mbox{and}\quad 
	\B_{\sigma}v_0=\B_{\sigma}(\A v_0)=0\quad\mbox{on}\ \ppp\Omega.
\end{cases}
\end{equation}
  Note that if we set $F(x,t)=f(t)R(x,t)$, then conditions \eqref{ft} 
and \eqref{Aeq1} are equivalent to \eqref{Aeq1(bis)}.
  Similarly, if we assume $p(x,t)=f(t)q(x,t)$, then conditions \eqref{ft} and 
\eqref{Aeq2} are equivalent to \eqref{Aeq2(bis)}.
  Now let us consider the following intermediate results.

\bigskip
\begin{lem}\label{l1}
  Let condition \eqref{Aeq1(bis)} be fulfilled.
  Then the IBVP \eqref{eq3} admits a unique solution 
$u\in\mathcal C([0,T];H^{2}(\Omega))$ satisfying
\[
	\A u\in\mathcal C([0,T];H^{2\gamma}(\Omega))
	\quad\mbox{and}\quad
	\ppp_t^\alpha u\in L^p(0,T;H^{2\gamma}(\Omega))
\]
  for all $0\le\gamma<1-1/(p\alpha)$.
  Moreover we have
\begin{align}\label{l1a}
	\|\A u\|_{\mathcal C([0,T];H^{2\gamma}(\Omega))}
		+\|\ppp_t^\alpha u\|_{L^p(0,T;H^{2\gamma}(\Omega))}
	\leq C\|F\|_{L^p(0,T;H^2(\Omega))}
\end{align}
  with $C>0$ depending on $\Omega$, $T$ and $\gamma$.
\end{lem}

\bigskip
\begin{lem}\label{l2}
  Let $F\in L^{\infty}(0,T;H^2(\Omega))$ satisfy $\B_{\sigma}F=0$ and 
condition 1) of \eqref{Aeq2(bis)} be fulfilled.
  Then the IBVP \eqref{eq4} admits a unique solution 
$v\in\mathcal C([0,T];H^{2}(\Omega))$ satisfying
\[
	\A v\in\mathcal C([0,T];H^{2\gamma}(\Omega))
	\quad\mbox{and}\quad
	\ppp_t^\alpha v\in L^{\infty}(0,T;H^{2\gamma}(\Omega))
\]
 for all $0\le\gamma<1$.
  Moreover, we have
\begin{equation}\label{l2a}
	\|\A v\|_{\mathcal C([0,T];H^{2\gamma}(\Omega))}
		+\|\ppp_t^\alpha v\|_{L^{\infty}(0,T;H^{2\gamma}(\Omega))}
	\leq C\|F\|_{L^{\infty}(0,T;H^2(\Omega))}
\end{equation}
  with $C$ depending on $\Omega$, $T$, 
$\|p\|_{L^{\infty}(0,T;H^{2}(\Omega))}$ and $\gamma$.
\end{lem}

\bigskip
\begin{lem}\label{l3}
  Let condition \eqref{Aeq2(bis)} be fulfilled.
  Then the IBVP \eqref{eq5} admits a unique solution 
$v\in\mathcal C([0,T];H^{2}(\Omega))$ satisfying
\[
	\A v\in\mathcal C([0,T];H^{2\gamma}(\Omega))
	\quad\mbox{and}\quad
	\ppp_t^\alpha v\in L^{\infty}(0,T;H^{2\gamma}(\Omega))
\]
  for all $0\le\gamma<1$.
  Moreover, we have
\begin{equation}\label{l3a}
	\|\A v\|_{\mathcal C([0,T];H^{2\gamma}(\Omega))}
	+\norm{\ppp_t^\alpha v}_{L^{\infty}(0,T;H^{2\gamma}(\Omega))}
	\leq C\|v_0\|_{H^4(\Omega)}
\end{equation}
  with $C$ depending on $\Omega$, $T$, 
$\|p\|_{L^{\infty}(0,T;H^{2}(\Omega))}$ and $\gamma$.
\end{lem}

\bigskip
  From these three lemmata we deduce easily Propositions \ref{p1} and \ref{p2}.

\subsection{Preliminary}
  We define the operator $A$ as $\A+1$ in $L^2(\Omega)$ equipped with the 
boundary condition $\B_{\sigma}h=0$;
\begin{equation}\label{A}
\begin{cases}
	D(A):=\{h\in H^2(\Omega);\ \B_{\sigma}h=0\ \mbox{on}\ \ppp\Omega \}, \\
	Ah:=\A h+h,\quad h\in D(A).
\end{cases}
\end{equation}
  Then $A$ is a selfadjoint and strictly positive operator with an orthonormal 
basis of eigenfunctions $(\phi_n)_{n\geq1}$ of finite order associated to an 
increasing sequence of eigenvalues $(\lambda_n)_{n\geq1}$.
  We can define the fractional power $A^{\gamma}$, $\gamma\ge0$, of $A$ by
\begin{equation}\label{frac}
\begin{aligned}
	&D(A^{\gamma})
	:=\left\{h\in L^2(\Omega);\ \sum_{n=1}^{\infty}\lambda_n^{2\gamma}
		\abs{(h,\phi_n)}^2<\infty\right\}, \\
	&A^{\gamma} h:=\sum_{n=1}^{\infty}\lambda_n^\gamma(h,\phi_n)\phi_n,
	\quad h\in D(A^{\gamma}).
\end{aligned}
\end{equation}
  Then $D(A^{\gamma})$ is a Hilbert space with the norm 
$\|\cdot\|_{D(A^{\gamma})}$ defined by 
$\|h\|_{D(A^{\gamma})}:=\|A^{\gamma}h\|$.
  Since $D(A)$ is continuously embedded into $H^2(\Omega)$ with norm 
equivalence (see Theorem 5.4 in Chapter 2 of \cite{LM} for example), we see by 
interpolation that
\begin{eqnarray*}
	&D(A^{\gamma})\subset H^{2\gamma}(\Omega), \\
	&C^{-1}\|h\|_{H^{2\gamma}(\Omega)}\le\|h\|_{D(A^{\gamma})}
		\le C\|h\|_{H^{2\gamma}(\Omega)},\quad h\in D(A^{\gamma})
\end{eqnarray*}
 for $0\le\gamma\le1$.

  In order to prepare for the arguments used in this paper, we consider the 
following Cauchy problem in $L^2(\Omega)$;
\begin{equation}\label{abs}
\begin{cases}
	\ppp_t^\alpha u(t)+A u(t)=F(t),\quad t\in(0,T), \\
	u(0)=0.
\end{cases}
\end{equation}
  We define the operator valued function $\{S_A(t)\}_{t\ge0}$ by
\[
	S_A(t)h
	=\sum_{n=1}^{\infty}(h,\phi_n)E_{\alpha,1}(-\lambda_n t^\alpha)\phi_n,
	\quad h\in L^2(\Omega),\ t\ge0,
\]
  with $E_{\alpha,\beta}$, $\alpha>0,\beta\in\R$, the Mittag-Leffler function 
given by 
\[
	E_{\alpha,\beta}(z)
	=\sum_{k=0}^{\infty}\frac{z^k}{\Gamma(\alpha k+\beta)}.
\]
  Recall that $S_A(t)\in W^{1,1}(0,T;\mathcal B(L^2(\Omega)))$ (e.g. \cite{BY} 
and \cite{SY}).
  Moreover, similarly to Theorem 2.2 in \cite{SY}, for 
$F\in L^{\infty}(0,T;L^2(\Omega))$, problem \eqref{abs} admits a unique 
solution given by 
\begin{equation}\label{sol}
	u(t)=\int_0^tA^{-1}S_A'(t-s)F(s)ds.
\end{equation}
  This solution is lying in $L^{\infty}(0,T;D(A^{\gamma}))$ for 
$0\le\gamma<1$, and, in view of Theorem 1 in \cite{BY}, we have 
\begin{equation}\label{es1}
	\|A^{\gamma-1}S_A'(t)h\|\le C t^{\alpha(1-\gamma)-1}\|h\|,
	\quad h\in L^2(\Omega),\ t>0.
\end{equation}
  In particular, the mapping $t\mapsto A^{-1}S_A'(t)$ belongs to 
$L^q(0,T;\B(L^2(\Omega)))$ if $q\in(1,\infty)$ satisfy 
$q(1-\alpha)<1$.
  Now we apply the following Young's inequality to \eqref{sol};

\bigskip
\begin{lem}\label{lem:yg}
  Let $f\in L^p(0,T)$ and $g\in L^q(0,T)$ with $1\le p,q\le\infty$ and 
$1/p+1/q=1$.
  Then the function $f*g$ defined by
\[
	f*g(t):=\int_0^tf(t-s)g(s)ds
\]
  belongs to $\mathcal C[0,T]$ and satisfies
\[
	|f*g(t)|\le\|f\|_{L^p(0,t)}\|g\|_{L^q(0,t)}, \quad t\in[0,T].
\]
\end{lem}

\bigskip
\begin{proof}
  Let $\tilde{f}$ and $\tilde{g}$ be defined by
\[
	\tilde{f}(t):=
\begin{cases}
	f(t),	&t\in(0,T), \\
	0,	&t\notin(0,T),
\end{cases}
	\quad\mbox{and}\quad
	\tilde{g}(t):=
\begin{cases}
	g(t),	&t\in(0,T), \\
	0,	&t\notin(0,T).
\end{cases}
\]
  Then applying the Young's inequality for functions on $\R$ (see Exercise 
4.30 in Brezis \cite{BR} or Appendix A in Stein \cite{ST} for example), we 
obtain the desired result.
\end{proof}

\bigskip
  Let $p\in(1,\infty]$ be as in \eqref{Aeq1(bis)}.
  Noting that $\A$ and $A^{-1}S_A'(t)$ commute, we see that for 
$F\in L^p(0,T;D(A))$,
\[
	\A u(t)=\int_0^tA^{-1}S_A'(t-s)\A F(s)ds.
\]
  By $p>1/\alpha$ and \eqref{es1}, the mapping $t\mapsto A^{-1}S_A'(t)$ 
belongs to $L^q(0,T;\mathcal B(L^2(\Omega)))$ where $q\in[1,\infty)$ satisfies 
$1/p+1/q=1$.
  Therefore by Lemma \ref{lem:yg}, $u$ belongs to 
$\mathcal C([0,T];D(A))$ and satisfies 
\begin{align}
	\|\A u(t)\|
&\le	\int_0^t\|A^{-1}S_A'(t-s)\|\|\A F(s)\|ds
\le 	C\int_0^t(t-s)^{\alpha-1}\|F(s)\|_{D(A)}ds \label{Au'} \\
&\le	C\left(\int_0^ts^{(\alpha-1)q}ds\right)^{1/q}
		\|F\|_{L^p(0,t;D(A))}
\le	Ct^{\alpha-1/p}\|F\|_{L^p(0,T;D(A))} \label{Au}.
\end{align}
  Thus we can define the map 
$\mathcal H:L^p(0,T;D(A))\to\mathcal C([0,T];D(A))$ by
\begin{equation}\label{Hwt}
	\mathcal H(w)(t):=\int_0^tA^{-1}S_A'(t-s)w(s)ds,
	\quad w\in L^p(0,T;D(A)).
\end{equation}
  By using these estimates, we will show the unique existence of the solution 
applying the fixed point theorem.

\subsection{Proof of Lemmata \ref{l1}-\ref{l3}}

\begin{proof}[\bf Proof of Lemma \ref{l1}]
  Let $A$ be the operator defined by \eqref{A}, then the IBVP 
\eqref{eq3} can be rewritten as 
\begin{equation}\label{l1c}
\begin{cases}
	\ppp_t^\alpha u(t)+Au(t)=u(t)+F(t),\quad t\in(0,T), \\
	u(0)=0,
\end{cases}
\end{equation}
  where $u(t):=u(\cdot,t)$ and $F(t):=F(\cdot,t)$.
  Noting that $F\in L^p(0,T;D(A))$ by \eqref{Aeq1(bis)}, we see from 
\eqref{sol} that the solution $u$ of \eqref{l1c} satisfies
\[
	u(t)=\mathcal H(u)(t)+\mathcal H(F)(t),\quad t\in(0,T),
\]
  where the map $\mathcal H$ is defined by \eqref{Hwt}.
  Therefore we will look for a fixed point of the map 
$\mathcal G:\mathcal C([0,T];D(A))\to\mathcal C([0,T];D(A))$ 
defined by
\begin{equation}\label{Gwt}
	\mathcal G(w)(t):=\mathcal H(w)(t)+\mathcal H(F)(t),
	\quad w\in\mathcal C([0,T]; D(A)),\ t\in(0,T).
\end{equation}
  By \eqref{Au'}, for $w\in\mathcal C([0,T]; D(A))$, we have
\begin{align*}
	\|\mathcal H(w)(t)\|_{D(A)}
&\le	C\int_0^t(t-s)^{\alpha-1}\|w(s)\|_{D(A)}ds
\le	C\left(\int_0^t(t-s)^{\alpha-1}ds\right)
		\|w\|_{\mathcal C([0,T];D(A))} \\
&=	\frac{Ct^{\alpha}}{\alpha}\|w\|_{\mathcal C([0,T];D(A))}.
\end{align*}
  Repeating the similar calculation, we get
\begin{align*}
	\|\mathcal H^2w(t)\|_{D(A)}
&\le	C\int_0^t(t-s)^{\alpha-1}\|\mathcal H w(s)\|_{D(A)}ds
\le	\frac{C}{\alpha}\left(\int_0^t(t-s)^{\alpha-1}s^{\alpha}ds\right)
		\|w\|_{\mathcal C([0,T];D(A))} \\
&=	\frac{C(\Gamma(\alpha)t^{\alpha})^2}{\Gamma(2\alpha+1)}
		\|w\|_{\mathcal C([0,T];D(A))}.
\end{align*}
  By induction, we have
\begin{equation}\label{Hnw}
	\|\mathcal H^nw(t)\|_{D(A)}
	\le\frac{C\left(\Gamma(\alpha)t^{\alpha}\right)^n}{\Gamma(n\alpha+1)}
		\|w\|_{\mathcal C([0,T];D(A))},
	\quad w\in\mathcal C([0,T];D(A)).
\end{equation}
  Therefore we obtain
\begin{align*}
	\|\mathcal G^n(w_1)-\mathcal G^n(w_2)\|
		_{\mathcal C([0,T];D(A))}
&=	\|\mathcal H^n(w_1-w_2)\|_{\mathcal C([0,T]; D(A))} \\
&\le	\frac{C\left(\Gamma(\alpha)T^{\alpha}\right)^n}{\Gamma(n\alpha+1)}
		\|w_1-w_2\|_{\mathcal C([0,T]; D(A))}
\end{align*}
  for $w_1,w_2\in\mathcal C([0,T]; D(A))$.
  Since $\mathcal G^n$ is a contraction for sufficiently large $n\in\N$, 
$\mathcal G$ admits a unique fixed point $u\in\mathcal C([0,T];D(A))
\subset\mathcal C([0,T];H^{2}(\Omega))$.
  Moreover we have
\[
	u=\mathcal G(u)=\mathcal G^n(u)
	=\mathcal H^n(u)+\sum_{k=1}^{n}\mathcal H^k(F)
\]
  for any $n\in\N$.
  Now we estimate each $\mathcal H^k(F)$.
  First, by \eqref{Au}, we have
\[
	\|\mathcal H(F)(t)\|_{D(A)}
\le	Ct^{\alpha-1/p}\|F\|_{L^p(0,T;D(A))}.
\]
  Next we apply \eqref{Au'} to have
\begin{align*}
	\|\mathcal H^2(F)(t)\|_{D(A)}
&\le	C\int_0^t(t-s)^{\alpha-1}\|\mathcal H(F)(s)\|_{D(A)}ds \\
&\le	C\left(\int_0^t(t-s)^{\alpha-1}s^{\alpha-1/p}ds\right)
		\|F\|_{L^p(0,T;D(A))} \\
&=	\frac{C\Gamma(\alpha)\Gamma(\alpha+1-1/p)}{\Gamma(2\alpha+1-1/p)}
		t^{2\alpha-1/p}\|F\|_{L^p(0,T;D(A))} \\
&\le	\frac{C\Gamma(\alpha)t^{2\alpha-1/p}}{\Gamma(2\alpha+1-1/p)}
		\|F\|_{L^p(0,T;D(A))}.
\end{align*}
  Repeating the similar calculation,
\begin{align*}
	\|\mathcal H^3(F)(t)\|_{D(A)}
&\le	C\int_0^t(t-s)^{\alpha-1}\|\mathcal H^2(F)(s)\|_{D(A)}ds \\
&\le	\frac{C\Gamma(\alpha)}{\Gamma(2\alpha+1-1/p)}
	\left(\int_0^t(t-s)^{\alpha-1}s^{2\alpha-1/p}ds\right)
		\|F\|_{L^p(0,T;D(A))} \\
&=	\frac{C\Gamma(\alpha)^2t^{3\alpha-1/p}}{\Gamma(3\alpha+1-1/p)}
		\|F\|_{L^p(0,T;D(A))}.
\end{align*}
  By induction, we obtain 
\[
	\|\mathcal H^k(F)\|_{\mathcal C([0,T];D(A))}
\le	\frac{C\Gamma(\alpha)^{k-1}T^{k\alpha-1/p}}{\Gamma(k\alpha+1-1/p)}
		\|F\|_{L^p(0,T;D(A))}.
\]
  Therefore
\begin{align*}
	\|u\|_{\mathcal C([0,T];D(A))}
&\le	\|\mathcal H^n(u)\|_{\mathcal C([0,T];D(A))}
	+\sum_{k=1}^{n}\|\mathcal H^k(F)\|_{\mathcal C([0,T];D(A))} \\
&\le	\frac{C\left(\Gamma(\alpha)T^{\alpha}\right)^n}{\Gamma(n\alpha+1)}
		\|u\|_{\mathcal C([0,T];D(A))}
		+\sum_{k=1}^{n}\frac{C\Gamma(\alpha)^{k-1}T^{k\alpha-1/p}}
		{\Gamma(k\alpha+1-1/p)}\|F\|_{L^p(0,T;D(A))}
\end{align*}
  and by taking sufficiently large $n\in\N$, we obtain
\begin{equation}\label{l1d}
	\|u\|_{\mathcal C([0,T]; D(A))}
	\le C\|F\|_{L^p(0,T;D(A))}
\end{equation}
  with $C$ depending on $T$ and $\Omega$.

  Now fix $0\le\gamma<1-1/(p\alpha)$.
  Then for all $t\in(0,T)$, we have $\A u(t)\in D(A^{\gamma})$ with
\[
	A^{\gamma}(\A u)(t)
=	\int_0^tA^{\gamma-1}S_{A}'(t-s)(\A u(s)+\A F(s))ds
\]
  and by \eqref{es1}, we have
\[
	\|A^{\gamma-1}S_{A}'(t)\|_{\B(L^2(\Omega))}\le Ct^{\mu-1},
\]
  where $\mu:=\alpha(1-\gamma)$.
  Since $\mu>1/p$, the mapping $t\mapsto A^{\gamma-1}S_{A}'(t)$ beongs to 
$L^q(0,T;\B(L^2(\Omega))$ where $q\in[1,\infty)$ satisfies $1/p+1/q=1$.
  Therefore $\A u$ belongs to $\mathcal C([0,T];D(A^{\gamma}))$ and 
\begin{align}\label{A(1+g)u}
	\|\A u(t)\|_{D(A^{\gamma})}
&=	\norm{\int_0^tA^{\gamma-1}S_{A}'(t-s)
		\left(\A u(s)+\A F(s)\right)ds} \nonumber \\
&\le	C\int_0^t(t-s)^{\mu-1}\|u(s)\|_{D(A)}ds
		+C\int_0^t(t-s)^{\mu-1}\|F(s)\|_{D(A)}ds \nonumber \\
&\le	C\left(\int_0^t(t-s)^{\mu-1}ds\right)
		\|u\|_{\mathcal C([0,T]; D(A))}
		+C\left(\int_0^t s^{q(\mu-1)}ds\right)^{1/q}
		\|F\|_{L^p(0,t;D(A))} \nonumber \\
&\le	CT^{\mu}\|u\|_{\mathcal C([0,T]; D(A))}
		+CT^{\mu-1/p}\|F\|_{L^p(0,T;D(A))}
\end{align}
  Combining this with \eqref{l1d}, we have
\[
	\|\A u(t)\|_{D(A^{\gamma})}
\le	C\|F\|_{L^p(0,T;D(A))}
\le	C\|F\|_{L^p(0,T;H^2(\Omega))}.
\]
  Hence we deduce that $\A u\in\mathcal C([0,T];H^{2\gamma}(\Omega))$ and
\[
	\|\A u\|_{\mathcal C([0,T];H^{2\gamma}(\Omega))}
\le	C\|F\|_{L^p(0,T;H^2(\Omega))}.
\]
  By the original equation $\ppp_t^\alpha u=-\A u+F$, we see that 
$\ppp_t^\alpha u$ belongs to $L^p(0,T;H^{2\gamma}(\Omega))$ with the estimate; 
\begin{align*}
	\|\ppp_t^\alpha u\|_{L^p(0,T;H^{2\gamma}(\Omega))}
&\le	C\|\A u\|_{L^p(0,T;H^{2\gamma}(\Omega))}
		+C\|F\|_{L^p(0,T;H^{2\gamma}(\Omega))} \\
&\le	C\|\A u\|_{\mathcal C([0,T];H^{2\gamma}(\Omega))}
		+C\|F\|_{L^p(0,T;H^{2}(\Omega))} \\
&\le	C\|F\|_{L^p(0,T;H^2(\Omega))},
\end{align*}
  which implies \eqref{l1a}.
  Thus we have completed the proof.
\end{proof}

\bigskip
  For the proof of Lemma \ref{l2}, we prepare the following fact;

\bigskip
\begin{lem}\label{pv}
  Let $u,v\in H^2(\Omega)$ and $d\le3$, then $uv\in H^2(\Omega)$ with the 
estimate
\[
	\|uv\|_{H^2(\Omega)}\le C\|v\|_{H^2(\Omega)}
\]
  with $C$ depending on $\|u\|_{H^2(\Omega)}$.
\end{lem}

\bigskip
  For this lemma, see Theorem 2.1 in Chapter II of Strichartz \cite{STR}.

\bigskip
\begin{proof}[\bf Proof of Lemma \ref{l2}]
  Similarly to Lemma \ref{l1}, the IBVP \eqref{eq4} can be rewritten as
\begin{equation}\label{eq6}
\begin{cases}
	\ppp_t^\alpha v(t)+Av(t)=(1-p(t))v(t)+F(t),\\
	v(0)=0,
\end{cases}
\end{equation}
  where $v(t):=v(\cdot,t)$ and $F(t):=F(\cdot,t)$.
  Moreover $p(t)$ denotes the multiplication operator by $p(x,t)$.
  Then we can see that the solution $v$ of \eqref{eq6} is a fixed point of the 
map $\mathcal K:\mathcal C([0,T];D(A))\to\mathcal C([0,T];D(A))$ 
defined by
\[
	\mathcal K(w)(t):=(\mathcal H(1-p(t))w)(t)+\mathcal H(F)(t),
	\quad w\in\mathcal C([0,T];D(A)),\ t\in(0,T).
\]
  Indeed, Lemma \ref{pv} and condition 1) of \eqref{Aeq2(bis)} yields that 
$(1-p)w$ belongs to $L^{\infty}(0,T;D(A))$ and satisfies
\[
	\|(1-p(t))w(t)\|_{D(A)}\le C\|w(t)\|_{D(A)}
\]
  with $C$ depending on $\|p\|_{L^{\infty}(0,T;H^2(\Omega))}$.
  Therefore we can see that $\mathcal K$ maps 
$\mathcal C([0,T];D(A))$.
  Moreover, by the similar calculation to \eqref{Hnw}, we have
\begin{equation}\label{H(1-p)}
	\|(\mathcal H(1-p))^n(w)\|_{\mathcal C([0,T]; D(A))}
	\le\frac{C(\Gamma(\alpha)t^{\alpha})^n}{\Gamma(n\alpha+1)}
		\|w\|_{\mathcal C([0,T]; D(A))},
	\quad w\in\mathcal C([0,T]; D(A))
\end{equation}
  and
\begin{equation}\label{H(1-p)Hk}
	\|(\mathcal H(1-p))^{n-1}(\mathcal H F)\|
		_{\mathcal C([0,T];D(A))}
	\le\frac{C(\Gamma(\alpha)t^{\alpha})^n}{\Gamma(n\alpha+1)}
		\|F\|_{L^{\infty}(0,T;D(A))},
	\quad F\in L^{\infty}(0,T;D(A)).
\end{equation}
  By \eqref{H(1-p)}, we find
\begin{align*}
	\|\mathcal K^n(w_1)-
		\mathcal K^n(w_2)\|_{\mathcal C([0,T];D(A))}
	\le\frac{C(\Gamma(\alpha)T^{\alpha})^n}{\Gamma(n\alpha+1)}
		\|w_1-w_2\|_{\mathcal C([0,T];D(A))}, \qquad\qquad\\
	w_1,w_2\in\mathcal C([0,T];D(A)),
\end{align*}
  which implies that $\mathcal K$ admits a unique fixed point 
$v\in\mathcal C([0,T];D(A))
\subset\mathcal C([0,T];H^{2}(\Omega))$.
  Then we have
\begin{equation}\label{vKv}
	v=\mathcal K(v)=\mathcal K^n(v)
	=(\mathcal H(1-p(t)))^n(v)
		+\sum_{k=1}^n(\mathcal H(1-p(t)))^{k-1}(\mathcal H F).
\end{equation}
  Repeating the argument in the proof of Lemma \ref{l1}, we deduce from 
\eqref{H(1-p)}, \eqref{H(1-p)Hk} and \eqref{vKv} that
\begin{equation}\label{l2b}
	\|v\|_{\mathcal C([0,T];D(A))}
	\le C\|F\|_{L^{\infty}(0,T;D(A))}
\end{equation}
  with $C$ depending on $T$, $\Omega$ and 
$\|p\|_{L^{\infty}(0,T;H^{2}(\Omega))}$.

  Next we fix $0\le\gamma<1$.
  Similarly to \eqref{A(1+g)u}, we have
\[
	\A v(t)\in D(A^{\gamma}),\quad t\in(0,T)
\]
  and
\begin{align*}
	\|\A v(t)\|_{D(A^{\gamma})}
&\le	C\norm{\int_0^tA^{\gamma-1}S_{A}'(t-s)
		\left((\A(1-p(s))v)(s)+\A F(s)\right)ds}  \\
&\le	C\int_0^t(t-s)^{\mu-1}\left(\|(1-p(s))v(s)\|_{D(A)}
		+\|F(s)\|_{D(A)}\right)ds \\
&\le	C\int_0^t(t-s)^{\mu-1}\left(\|v(s)\|_{D(A)}
		+\|F(s)\|_{D(A)}\right)ds
\end{align*}
  with $\mu=\alpha(1-\gamma)$.
  Therefore $\A v$ belongs to $\mathcal C([0,T];H^{2\gamma}(\Omega))$ and 
satisfies
\begin{align*}
	\|\A v\|_{\mathcal C([0,T]; H^{2\gamma}(\Omega))}
&\le	\|\A v\|_{\mathcal C([0,T];D(A^{\gamma}))}
\le	CT^{\mu}\left(\|v\|_{\mathcal C([0,T];D(A))}
		+\|F\|_{L^{\infty}(0,T;D(A))}\right) \\
&\le	C\|F\|_{L^{\infty}(0,T;D(A))}
\le	C\|F\|_{L^{\infty}(0,T;H^2(\Omega))},
\end{align*}
  where we have used \eqref{l2b}.
  Moreover, combining this with the original equation, we also have 
$\ppp_t^\alpha v\in L^{\infty}(0,T;H^{2\gamma}(\Omega))$ and \eqref{l2a}.
\end{proof}

\bigskip
\begin{proof}[\bf Proof of Lemma \ref{l3}]
  We split the solution $v$ of \eqref{eq5} into two terms $v=w+v_0$ where $w$ 
solves
\begin{equation}\label{eq7}
\left\{\begin{aligned}
	\ppp_t^\alpha w(x,t)+\A w(x,t)+p(x,t)w(x,t)=F(x,t),\quad &(x,t)\in Q,\\
	\B_{\sigma} w(x,t)=0,\quad &(x,t)\in\Sigma, \\
	w(x,0)=0,\quad &x\in\Omega
\end{aligned}\right.
\end{equation}
  with $F(x,t):=-(\A+p(x,t))v_0(x)$.
  Then \eqref{Aeq2(bis)} implies $F\in L^{\infty}(0,T;D(A))$ with 
the estimate
\[
	\|F\|_{L^{\infty}(0,T;H^2(\Omega))}\leq C\|v_0\|_{H^4(\Omega)}.
\]
  By Lemma \ref{l2}, the IBVP \eqref{eq7} admits a unique solution 
$w\in \mathcal C([0,T];H^{2}(\Omega))$ satisfying
\[
	\A w\in\mathcal C([0,T];H^{2\gamma}(\Omega))
	\quad\mbox{and}\quad
	\ppp_t^\alpha w\in L^{\infty}(0,T;H^{2\gamma}(\Omega)).
\]
  Moreover
\[
	\|\A w\|_{\mathcal C([0,T];H^{2\gamma}(\Omega))}
		+\|\ppp_t^\alpha w\|_{L^{\infty}(0,T;H^{2\gamma}(\Omega))}
	\leq C\|F\|_{L^{\infty}(0,T;H^2(\Omega))}
	\leq C\|v_0\|_{H^4(\Omega)}.
\]
  Therefore the IBVP \eqref{eq5} admits a unique solution 
$v\in\mathcal C([0,T];H^{2}(\Omega))$ satisfying
\[
	\A v\in\mathcal C([0,T];H^{2\gamma}(\Omega))
	\quad\mbox{and}\quad
	\ppp_t^\alpha v\in L^{\infty}(0,T;H^{2\gamma}(\Omega)).
\]
  From the above estimate, we deduce \eqref{l3a}.
\end{proof}

\section{Proof of Theorem \ref{t1}}

  In this section, we prove Theorem \ref{t1}.
  To this end, we prepare the following lemmata with Gronwall type 
inequalities;

\bigskip
\begin{lem}\label{lem:h1}
  Let $C,\alpha>0$ and $u,d\in L^1(0,T)$ be nonnegative functions satisfying
\[
	u(t)\le Cd(t)+C\int_0^t(t-s)^{\alpha-1}u(s)ds,\quad t\in(0,T),
\]
  then we have
\[
	u(t)\le Cd(t)+C\int_0^t(t-s)^{\alpha-1}d(s)ds,\quad t\in(0,T).
\]
\end{lem}

\bigskip
  For the proof, see Lemma 7.1.1 p.188 of \cite{H}.

\bigskip
\begin{lem}\label{gron}
  We take $2\le p\le\infty$ and $\mu>2/p$.
  Let $f\in L^{\infty}(0,T)$ and $u,R\in L^p(0,T)$ be non-negative functions 
satisfying the integral inequality
\begin{equation}\label{eq:gron}
	f(t)\le u(t)+\int_0^t(t-s)^{\mu-1}f(s)R(s)ds,
	\quad\mbox{a.e.}\ t\in(0,T).
\end{equation}
  Then we have
\begin{equation}\label{fu}
	\|f\|_{L^p(0,T)}\le C\|u\|_{L^p(0,T)},
\end{equation}
  where the constant $C$ depends on $p$, $\mu$, $T$ and $\|R\|_{L^p(0,T)}$.
\end{lem}

\bigskip
\begin{proof}
  We set $d(t):=\|f\|_{L^p(0,t)}^p$.
  From equation \eqref{eq:gron}, we have
\[
	|f(s)|^p
	\le C|u(s)|^p+C\left|\int_0^s(s-\xi)^{\mu-1}f(\xi)R(\xi)d\xi\right|^p,
\]
  which implies
\begin{equation}\label{eq:dt}
	d(t)\le C\|u\|_{L^p(0,T)}^p
	+C\int_0^t\left|\int_0^s(s-\xi)^{\mu-1}f(\xi)R(\xi)d\xi\right|^pds.
\end{equation}
  Now we estimate the right-hand side of the above.
  By the Cauchy-Schwarz inequality, 
\begin{align*}
	\int_0^s|f(\xi)R(\xi)|^{p/2}d\xi
=	\int_0^s|f(\xi)|^{p/2}\cdot|R(\xi)|^{p/2}d\xi
\le	\left(\int_0^s|f(\xi)|^pd\xi\right)^{1/2}
		\left(\int_0^s|R(\xi)|^pd\xi\right)^{1/2},
\end{align*}
  that is,
\[
	\|fR\|_{L^{p/2}(0,s)}
\le	\|f\|_{L^p(0,s)}\|R\|_{L^p(0,s)}.
\]
  Therefore if $p>2$, then Lemma \ref{lem:yg} yields that
\begin{align*}
	\left|\int_0^s(s-\xi)^{\mu-1}f(\xi)R(\xi)d\xi\right|
\le	\left(\int_0^s\xi^{r(\mu-1)}ds\right)^{1/r}\|fR\|_{L^{p/2}(0,s)}
\le	C\|f\|_{L^p(0,s)}\|R\|_{L^p(0,s)},
\end{align*}
  where $r\in[1,\infty)$ satisfies $2/p+1/r=1$.
  For $p=2$, we also have
\[
	\left|\int_0^s(s-\xi)^{\mu-1}f(\xi)R(\xi)d\xi\right|
\le	s^{\mu-1}\|fR\|_{L^1(0,s)}
\le	C\|f\|_{L^2(0,s)}\|R\|_{L^2(0,s)}.
\]
  Thus for any $p\ge2$, we have
\begin{equation}\label{eq:fp}
	\left|\int_0^s(s-\xi)^{\mu-1}f(\xi)R(\xi)d\xi\right|\le Cd(s),
\end{equation}
  where $C$ depends on $T$, $p$, $\mu$ and $\|R\|_{L^p(0,T)}$.
  By \eqref{eq:dt} and \eqref{eq:fp}, we have
\[
	d(t)\le C\|u\|_{L^p(0,T)}^p+C\int_0^td(s)ds,\quad t\in(0,T).
\]
  Hence by the Gronwall inequality, we have
\[
	d(t)\le C\|u\|_{L^p(0,T)}^p,\quad t\in(0,T)
\]
  with $C$ depending on $p$, $\mu$, $T$ and $\|R\|_{L^p(0,T)}$.
  Thus we have proved \eqref{fu}.
\end{proof}

\bigskip
  Now we are ready to prove Theorem \ref{t1}.

\bigskip
\begin{proof}[\bf Proof of Theorem \ref{t1}]
  Let $u_i$ be the solutions to \eqref{eq1} corresponding to $f_i$ ($i=1,2$) 
and set $u:=u_1-u_2$ and $f:=f_1-f_2$.
  Then $u$ solves \eqref{eq1} and is given by
\[
	u(t)=\int_0^tA^{-1}S_{A}'(t-s)u(s)+\int_0^tA^{-1}S_{A}'(t-s)f(s)R(s)ds,
\]
  where $u(t):=u(\cdot,t)$ and $R(t):=R(\cdot,t)$.

  First we estimate $\|u(t)\|_{D(A)}$.
  Similarly to the calculation in \eqref{Au'}, we have
\begin{align}\label{t1c}
	\|u(t)\|_{D(A)}
&\le	C\int_0^t(t-s)^{\alpha-1}\|u(s)\|_{D(A)}ds
	+C\int_0^t(t-s)^{\alpha-1}|f(s)|\|R(s)\|_{D(A)}ds \nonumber \\
&=	C\int_0^t(t-s)^{\alpha-1}\|u(s)\|_{D(A)}ds+Cd(t),
\end{align}
  where we have set 
\[
	d(t):=
	\int^t_0(t-s)^{\alpha-1}|f(s)|\|R(s)\|_{D(A)}ds.
\]
  Applying Lemma \ref{lem:h1} to \eqref{t1c}, we have
\begin{align}\label{t1c'}
	\|u(t)\|_{D(A)}
\le 	Cd(t)+C\int^t_0 (t-s)^{\alpha-1} d(s) ds,\quad 0<t<T.
\end{align}
  Here for $\nu>0$, we note 
\begin{align}\label{t1e}
	\int^t_0 (t-s)^{\nu-1}d(s)ds
&=	\int^t_0 (t-s)^{\nu-1}\left(\int^s_0(s-\xi)^{\alpha-1}|f(\xi)|
		\|R(\xi)\|_{D(A)}d\xi\right)ds \nonumber \\
&=	\int^t_0 \left( \int^t_{\xi}(t-s)^{\nu-1}(s-\xi)^{\alpha-1}
		ds\right)|f(\xi)|\|R(\xi)\|_{D(A)}d\xi \nonumber \\
&=	B(\nu,\alpha)\int^t_0 (t-\xi)^{\nu+\alpha-1}
		|f(\xi)|\|R(\xi)\|_{D(A)}d\xi
\end{align}
  where $B(\cdot,\cdot)$ is the Beta function.
  In particular, for $\nu=\alpha$, we have
\begin{align*}
	\int^t_0 (t-s)^{\alpha-1}d(s) ds
&= 	B(\alpha,\alpha)\int^t_0(t-s)^{2\alpha-1}|f(s)|
		\|R(s)\|_{D(A)}ds \\
&\le 	T^{\alpha}B(\alpha,\alpha)
		\int^t_0(t-s)^{\alpha-1}|f(s)|\|R(s)\|_{D(A)}ds \\
&\le	Cd(t).
\end{align*}
  Hence the following estimate follows from \eqref{t1c'};
\begin{equation}\label{t1f}
	\|u(t)\|_{D(A)}\le Cd(t),\quad 0<t<T.
\end{equation}

  Next we estimate $\|\A u(t)\|_{D(A^{\gamma})}$ for 
$d/4<\gamma<1-2/(p\alpha)$.
  Repeating the calculation in \eqref{A(1+g)u}, we find
\begin{align*}
	\|\A u(t)\|_{D(A^{\gamma})}
\le	C\int_0^t(t-s)^{\mu-1}\left(\|u(s)\|_{D(A)}
		+|f(s)|\|R(s)\|_{D(A)}\right)ds,
	\quad\mbox{a.e.}\ t\in(0,T)
\end{align*}
  where $\mu=\alpha(1-\gamma)$.
  By \eqref{t1e} with $\nu=\mu$ and \eqref{t1f}, we obtain
\begin{align*}
	\|\A u(t)\|_{D(A^{\gamma})}
&\le	C\int^t_0 (t-s)^{\mu-1}d(s)ds
		+ C\int^t_0 (t-s)^{\mu-1}|f(s)|\|R(s)\|_{D(A)}ds \\
&=	CB(\mu,\alpha)\int^t_0 (t-s)^{\mu+\alpha-1}
		|f(s)|\|R(s)\|_{D(A)}ds \nonumber \\
&\qquad	+C\int^t_0 (t-s)^{\mu-1}|f(s)|\|R(s)\|_{D(A)}ds \\
&\le	CT^{\alpha}B(\mu,\alpha)\int^t_0(t-s)^{\mu-1}
		|f(s)|\|R(s)\|_{D(A)}ds \\
&\qquad	+ C\int^t_0 (t-s)^{\mu-1}|f(s)|\|R(s)\|_{D(A)}ds  \\
&\le	C\int^t_0 (t-s)^{\mu-1}|f(s)|\|R(s)\|_{D(A)}ds.
\end{align*}

  Finally we estimate $|\A u(x_0,t)|$ and complete the proof.
  Since $\gamma>d/4$, the Sobolev embedding theorem yields
\begin{align}\label{t1g}
	|\A u(x_0,t)|
\le	C\|\A u(\cdot,t)\|_{H^{2\gamma}(\Omega)}
\le	C\|\A u(t)\|_{D(A^{\gamma})}
\le	C\int^t_0 (t-s)^{\mu-1}|f(s)|\|R(s)\|_{D(A)}ds.
\end{align}
  From the original equation, we get
\begin{equation}\label{ori}
	f(t)R(x_0,t)=\ppp_t^{\alpha}u(x_0,t)+\A u(x_0,t),
	\quad\mbox{a.e.}\ t\in(0,T).
\end{equation}
  Combining this with \eqref{t1a} and \eqref{t1g}, we get
\begin{align}\label{f}
	|f(t)|
&\le 	\frac{1}{\delta}\left(|\ppp_t^{\alpha}u(x_0,t)|+|\A u(x_0,t)|\right) 
		\nonumber \\
&\le 	C|\ppp_t^{\alpha}u(x_0,t)|
		+C\int^t_0 (t-s)^{\mu-1}|f(s)|\|R(s)\|_{D(A)}ds,
	\quad\mbox{a.e.}\ t\in(0,T)
\end{align}
  with $C$ depending on $\delta$, $\Omega$ and $T$.
  By Lemma \ref{gron}, we see that
\[
	\|f\|_{L^p(0,T)}\le C\|\ppp_t^{\alpha}u(x_0,\cdot)\|_{L^p(0,T)},
\]
  which implies \eqref{t1b}.
  Moreover, by \eqref{t1g} and \eqref{ori}, we have
\begin{align*}
	|\ppp_t^{\alpha}u(x_0,t)|
&\le	|f(t)R(x_0,t)|+|\A u(x_0,t)| \\
&\le	C|f(t)|\|R(\cdot,t)\|_{H^{2}(\Omega)}
		+C\int^t_0 (t-s)^{\mu-1}|f(s)|\|R(s)\|_{D(A)}ds \\
&\le	C\|f\|_{L^{\infty}(0,T)}\|R(t)\|_{D(A)}
		+C\|f\|_{L^{\infty}(0,T)}\int^t_0 (t-s)^{\mu-1}
		\|R(s)\|_{D(A)}ds \\
&\le	C\|f\|_{L^{\infty}(0,T)}\|R(t)\|_{D(A)}
	+C\|f\|_{L^{\infty}(0,T)}T^{\mu-1/p}\|R\|_{L^p(0,T;D(A))},
	\quad\mbox{a.e.}\ t\in(0,T).
\end{align*}
  Therefore,
\begin{align*}
	\|\ppp_t^{\alpha}u(x_0,\cdot)\|_{L^p(0,T)}
&\le	C\|f\|_{L^{\infty}(0,T)}\|R\|_{L^p(0,T;D(A))}
	+C\|f\|_{L^{\infty}(0,T)}T^{\mu}\|R\|_{L^p(0,T;D(A))} \\
&\le	C\|f\|_{L^{\infty}(0,T)}.
\end{align*}
  Thus we have proved \eqref{t1b'}.
\end{proof}

\section{Proof of Theorem \ref{t2}}
  In this section, we prove Theorem \ref{t2}.
  We first prepare the following generalized Gronwall's inequality;

\bigskip
\begin{lem}\label{lem:h2}
  Let $\mu,a,b>0$ and $f\in L^1(0,T)$ be nonnegative function satisfying 
the integral inequality
\[
	f(t)\le a+b\int_0^t(t-s)^{\mu-1}f(s)ds, \quad\mbox{a.e.}\ t\in(0,T).
\]
  Then we have
\[
	f(t)\le aE_{\mu,1}\left((b\Gamma(\mu))^{1/\mu}t^{\mu}\right),
	\quad\mbox{a.e.}\ t\in(0,T).
\]
\end{lem}

\bigskip
  For the proof, see Lemma 7.1.2 on p.189 of \cite{H}.
  Now we are ready to prove Theorem \ref{t2}.

\bigskip
\begin{proof}[\bf Proof of Theorem \ref{t2}]
  Let $v_i$ be the solutions to \eqref{eq2} corresponding to $f_i$ ($i=1,2$) 
and set $v:=v_1-v_2$ and $f:=f_2-f_1$.
  Then $v$ solves \eqref{eq4} with $p(x,t)=f_1(t)q(x,t)$ and 
$F(x,t)=f(t)q(x,t)v_2(x,t)$.
  Recall that $v$ is given by
\[
	v(t)=\int_0^t A^{-1}S_{A}'(t-s)((1-p(t))v)(s)
		+\int_0^tf(s)A^{-1}S_{A}'(t-s)R(s)ds,
\]
  where we have set $v(t):=v(\cdot,t)$ and $R(t):=q(\cdot,t)v_2(\cdot,t)$.
  Moreover, $p(t)$ denotes the multiplication operator by 
$p(x,t):=f_1(t)q(x,t)$.

  First we estimate $\|v(t)\|_{D(A)}$.
  Since $(1-p(t))v(t),R(t)\in D(A)$ by \eqref{Aeq2}, we repeat the 
calculation in \eqref{A(1+g)u} to have
\begin{align*}\label{t2c}
	\|v(t)\|_{D(A)}
&\le 	C\int^t_0 (t-s)^{\alpha-1}\|(1-p(t))v(s)\|_{D(A)}ds
		+C\int^t_0 (t-s)^{\alpha-1}|f(s)|\|R(s)\|_{D(A)}ds 
		\nonumber \\
&\le 	C\int^t_0 (t-s)^{\alpha-1}\|v(s)\|_{D(A)}ds
		+C\int^t_0 (t-s)^{\alpha-1}|f(s)|ds.
\end{align*}
  with $C$ depending on $\Omega$, $M$ and 
$\|q\|_{L^{\infty}(0,T;H^{2}(\Omega))}$.
  Then repeating the arguments used in Theorem \ref{t1}, we obtain
\[
	\|v(t)\|_{D(A)}
\le 	C\int^t_0 (t-s)^{\alpha-1}|f(s)|ds,\quad 0<t<T.
\]
  and from this estimate we also deduce that for any $0\le\gamma<1$,
\[
	\|\A v(t)\|_{D(A^{\gamma})}
\le	C\int^t_0 (t-s)^{\mu-1}|f(s)|ds,\quad 0<t<T,
\]
  where $\mu:=\alpha(1-\gamma)$.
  Therefore by taking $\gamma\in(d/4,1)$, we have
\begin{align}\label{t2g}
	|\A v(x_0,t)+p(x_0,t)v(x_0,t)|
&\le	C\|\A v(\cdot,t)+p(\cdot,t)v(\cdot,t)\|_{H^{2\gamma}(\Omega)} 
		 \nonumber \\
&\le	C\|\A v(\cdot,t)\|_{H^{2\gamma}(\Omega)}
		+C\|v(\cdot,t)\|_{H^{2\gamma}(\Omega)} \nonumber \\
&\le	C\|\A v(t)\|_{D(A^{\gamma})}
\le	C\int^t_0 (t-s)^{\mu-1}|f(s)|ds.
\end{align}
  From the original equation, we have
\begin{equation}\label{fR}
	f(t)R(x_0,t)=
	\ppp_t^{\alpha}v(x_0,t)+\A v(x_0,t)+p(x_0,t)v(x_0,t),
	\quad\mbox{a.e.}\ t\in(0,T).
\end{equation}
  On the other hand, from \eqref{t2a}, we deduce that
\[
	|R(x_0,t)|\geq c>0,\quad\mbox{a.e.}\ t\in(0,T)
\]
  with $c$ depending on $\delta$, $\Omega$ and $T$.
  Therefore, combining this with \eqref{t2g} and \eqref{fR}, we obtain
\begin{align*}
	|f(t)|
&\le	C|\ppp_t^{\alpha}v(x_0,t)|+C|\A v(x_0,t)+p(x_0,t)v(x_0,t)|\\
&\le	C\|\ppp_t^{\alpha}v(x_0,\cdot)\|_{L^{\infty}(0,T)}
		+C\int^t_0 (t-s)^{\mu-1}|f(s)|ds,
	\quad\mbox{a.e.}\ t\in(0,T).
\end{align*}
  Applying Lemma \ref{lem:h2}, we see that
\[
	|f(t)|\le C\|\ppp_t^{\alpha}v(x_0,\cdot)\|_{L^{\infty}(0,T)}.
\]
  Thus we have proved the second inequality in \eqref{t2b}.
  Moreover, by \eqref{fR}, we have
\begin{align*}
	|\ppp_t^{\alpha}v(x_0,t)|
&\le 	|f(t)R(x_0,t)|+|\A v(x_0,t)+p(x_0,t)v(x_0,t)| \\
&\le	|f(t)|\|R(\cdot,t)\|_{D(A)}+C\int^t_0 (t-s)^{\mu-1}|f(s)|ds \\
&\le	C\left(\|R\|_{L^{\infty}(0,T;D(A))}
		+\frac{T^{\mu}}{\mu}\right)\|f\|_{L^{\infty}(0,T)}.
\end{align*}
  Thus we have proved the first inequality in \eqref{t2b}.
\end{proof}

\section*{Acknowledgements}

  Both authors would like to express thank to Professor Masahiro Yamamoto, who 
is the supervisor of the first named author, for his useful comments and other 
kind helps.
  The first named author is granted by the FMSP program at Graduate School of 
Mathematical Sciences of The University of Tokyo.

\end{document}